\documentclass{amsart}
\usepackage{amssymb, mathtools, fullpage, color, cite, url}
\usepackage[colorlinks=true, pdfstartview=FitV,linkcolor=blue,citecolor=blue,urlcolor=blue]{hyperref}
\usepackage[a4paper]{geometry}

\begin{document}

\title{Mathematics: the Rise of the Machines}

\author[Y.-H. He]{Yang-Hui He}
\address{London Institute for Mathematical Sciences, Royal Institution, London W1S 4BS, UK \&
Merton College, University of Oxford, OX14JD, UK}
\email{\href{mailto:hey@maths.ox.ac.uk}{hey@maths.ox.ac.uk}}

\begin{abstract}
We argue how AI can assist mathematics in three ways: theorem-proving, conjecture formulation, and language processing.
Inspired by initial experiments in geometry and theoretical physics in 2017, we summarize how this emerging field has grown over the past years, and show how various machine-learning algorithms can help with pattern detection across disciplines in the mathematical sciences. 
At the heart is the question how does AI help with theoretical discovery, and the implications for the future of mathematics.
~\\

Invited contribution to the {\it Nieuw Archief voor Wiskunde} of the {\it Koninklijk Wiskundig Genootschap} and based on a recent public lecture of the author to the Royal Institution of Great Britain \cite{riYHH}. 
\end{abstract}
    
\date{November 2025}

\maketitle

Whether you are of the opinion that Artificial Intelligence (AI) is a threat to humanity or of the opposite belief that it is mere hype, it cannot be denied that - whatever its precise definition - AI is already transforming our daily lives in almost every respect.
Even though AI has a much older history than many would imagine, the formidable combination of compute (your mobile phone out-powering the entire lunar programme in the 1960s) and the internet (essentially all of human knowledge at your finger-tips), has galvanized this transformative impact in the last two decades.
The word ``Artificial Intelligence'' was first coined \footnote{Interestingly, the initial wording was ``Computational Intelligence'', and later changed to AI because it sounded more impressive, as grant proposals should.} in a grant proposal in 1956 by McCarthy et al.~\cite{mccarthy2006proposal}.
Shortly thereafter, the ``Perceptron'' - the first artificial neural network (NN) - was established, as an impressive array of photo-receptors and wires \cite{rosenblatt1957perceptron}.

While I shall leave the question of the influence of AI on our civilization to the experts in AI research, to sociologists, and to philosophers, it is my goal here to give a rapid overview of the emergent but blossoming field of {\em AI for Mathematics} (AI4Maths). 
To be fair, both ``mathematics and theoretical physics for AI'' and ``AI for the experimental sciences'' have long and distinguished histories.
For instance, early universal approximation theorems \cite{cybenko1989approximation}
 that ensure good estimations from NNs, or the Boltzman machine as an energy-based model for cognition \cite{ackley1985learning}, date to the late 1980s.
 So too, did CERN's first usage of machine-learning (ML) for detecting fundamental particles \cite{perret1990new}. 
The award of the 2024 Nobel Prize in Physics to Hinton and Hopfield, and in Chemistry to Hassabis and Jumper, is proof of the establishment's timely recognition of the r\^{o}le of AI in scientific discovery.

Pure mathematics and theoretical physics, by their very nature of rigour, may seem unlikely companions to the big data and statistical inferences inherent to AI.
Yet, I will describe here how in three ways AI will be a positive driving force for mathematical research: bottom-up, top-down, and meta-mathematics.
Admittedly, these three words have various precise meanings in the philosophy of mathematics and in the philosophy of science, I beg the reader's indulgence for borrowing them loosely and shall define them in the ensuing.
The interested reader is referred to some recent and more in-depth reviews \cite{He:2021oav,douglas2022machine,castelvecchi2023will,Gukov:2024buj,he2024can,he2024ai,liang2024mathematics,douglas2025mathematical,kontorovich2025shape} which attempt to describe the idea of {\it AI4Maths}.

\section{Bottom-Up Mathematics: }
Typically, one thinks of mathematics as an austere three-step process: (1) laying out the {\it axioms} and {\it definitions} of terms and concepts with precision; (2) claiming a statement involving the terms; and (3) using logic to deduce the statement.
Since Euclid's immortal tome {\it The Elements}, this is how mathematics is always presented. A classic example (Proposition 20 of Book IX {\it cit.~ibid.}) is as follows.
(1) Definition: a prime number is a positive integer which is only divisible by 1 and itself; (2) Proposition: there is an infinite number of such prime numbers; (3) Proof: Standard \footnote{Suppose there is a finite number $n$ of primes, let them be ordered as $p_1 = 2, p_2 =3, \ldots, p_n$. Consider the product $N := p_1 \cdot p_2 \ldots \cdot p_n$ and the number $N+1$. Now, $N+1$ is either prime or not. If it is a prime, then we have found one larger than $p_n$, contradicting our assumption that $p_n$ is the largest prime. Thus $N+1$ must be composite, whereby it must be divisible by some prime $q$ between $p_1$ and $p_n$. But $N$ itself is divisible by $q$ by construction, being the product over all primes. Hence their difference - which is equal to 1 - must be divisible by $q \geq 2$, another contradiction. Hence reductio ad absurdum implies our initial assumption that $p_n$ exists is false [QED]. Of course, there are some further properties in the arithmetic of integers which are used here; they need to be and can be proven separately.}.

This ``bottom-up'' presentation of statements in mathematics, building sentence by sentence from the foundational definitions and axioms, perhaps reached its zenith with the publication \cite{PM} of {\it Principia Mathematica} in 1910.
With a sparsity of words, a plethora of symbols, and an inhuman coldness, the monumental work was famously able to reach ``1+1=2'' on page 362.
Even though this entire programme of {\em formalizing mathematics} was in part killed by G\"odel-Church-Turing \cite{godel1931formal,church1936unsolvable,turing1936computable} with their showing the existence of undecidable and uncomputable statements, researchers were not deterred in using the earliest version of the computer - the Logical Theory Machine - to formalize and prove a number of propositions of the {\it Principia} \cite{newell1957empirical}. This was almost immediately after the appearance of the first NN \cite{rosenblatt1957perceptron}.
Likewise, the idea of ``mechanized mathematics'' \cite{wang1960toward} dates to that time.

Today's answer to this tradition of ``bottom-up'' mathematics is the automated proof-assistant \cite{Tao_AI}, notably {\em Mathlib}, built in the {\em Lean} programming language \cite{The_mathlib_Community_2020,de2015lean}.
Launched in 2013, Project Xena \cite{xena} has completed the formalization in {\em Lean}, of the  statements and proofs of essentially all undergraduate-level mathematics. 
Equally impressive is the fact that formalizing whilst doing collaborations has already led to major results by the top mathematicians
\cite{scholze2022liquid,gowers2023conjecturemarton}.

As proof co-pilots are quickly becoming the norm for mathematical research, it should be emphasized that so far there is no large-scale AI involved in {\em Mathlib}.
The main challenge here is one of engineering: there is only on the order of millions of lines of {\em Lean}, whereas currently reliable language models (to which I shall turn later) typically require billions of lines to train. It is conceivable that in the not too distant future we will reach auto-formalization which will take all journal papers in LaTeX and generate verifiable {\it Lean} code for the whole corpus.
Already, Deepmind's AlphaGeo2 \cite{chervonyi2025goldmedalistperformancesolvingolympiad} and AlphaProof \cite{castelvecchi2024deepmind} have {\it Lean} incorporated into their core algorithms.
Once accurate auto-formalization is achieved, it might no longer be a fantasy that proof-paths maybe found, using reinforcement learning, for major open problems.

\section{Top-Down Mathematics: }
Lest the reader be led to think that mathematical research is always done in this formal and rather stark ``bottom-up'' way, I must point out that in practice, the opposite is the case.
To quote the great V.~Arnol'd, ``mathematics is a branch of physics where the experiments are cheap.''
Indeed, even rigorous proofs are formulated first by some ``gut-feeling'', as summarized in this wonderful article by Thurston \cite{thurston2006proof}.
Certainly, theoretical physics, the closest cousin to pure mathematics, is done this way. As a famous example, we have achieved 19 significant-digit agreement between theory and experiment in the standard model of particle physics, but we still do not have a precise definition of the underlying quantum field theory.
This fearless bravado is the very reason why string theory has produced so many new insights in pure mathematics.

The idea of {\em experimental mathematics} is as old as logic and reasoning. It is how the practitioner actually does theoretical research: we play around, trying to get intuition, and then formulate the right definitions and proofs, before finally writing up the papers more formally.
In today's parlance, this might be called {\em data-driven} mathematics, where one builds experience by finding patterns and by a broad understand of the research landscape, looking at the subject from ``top-down''.
After all, in the beautiful words of G.~H.~Hardy, ``A mathematician, like a painter or poet, is a maker of patterns. If his patterns are more permanent than theirs, it is because they are made with ideas.'' Mathematical data is distinctly different from data from the the natural world, it is without statistical error, and often represented by integers rather than reals whence also has no freedom for numerical error. These have been dubbed {\it Platonic Data} \cite{douglas2025mathematical}, whose machine-learning can uncover underlying mathematical structure \cite{He:2021oav}.

Continuing with our example of prime numbers, one of the high-lights of this data-driven mathematics was the 16 year-old Gauss defining the prime-counting function $\pi(x) := \#\{p \mbox{ prime }: p \leq (x \in  \mathbb{R}_+)\}$, the number of primes not exceeding a real positive $x$. He was able to eyeball \footnote{Actually, he invented regression and the earliest notions of statistics in order to gain this intuition.} that $\pi(x) \sim x / \log(x)$.
The proof had to wait some 50 years, especially for the development of complex analysis, before it was established as the Prime Number Theorem, one of the cornerstones of mathematics. Based on experience and intuition, Gauss {\em formulated a conjecture} by looking at Platonic data, whereby establishing a good problem.
In some sense, ``bottom-up'' mathematics can be thought of as finding proofs, and ``top-down'' mathematics, as finding interesting problems.
One needs to bear in mind that two of the remaining six Millennium Problems \cite{carlson2023millennium} - the Riemann Hypothesis and the Birch--Swinnerton-Dyer (BSD) Conjecture - arose from mathematical experimentation by listing the first zeroes or by plotting rank and conductor.

Starting from the late 1980s, attempts were made to automate the process of formulating conjectures in combinatorics \cite{davila2024automated}, which in 2017 took the updated form as {\em TxGraffiti}.
In the same year, the first AI-assisted mathematical experiments were performed \cite{He:2017aed}, where datasets from algebraic geometry were passed through a neural-network, in order to uncover potential new structures. Therein, calculations in algebraic geometry - and later \cite{He:2018jtw}, all mathematical calculations in one way or another -  were interpreted as {\em image processing}.
So too, did 2017 see machine-learning applied to mathematical physics \cite{Carifio:2017bov,Krefl:2017yox,Ruehle:2017mzq} in exploring the string theory landscape \footnote{As indeed, did Sophia become the first robot to receive a human passport \cite{weller2017meet}.}.
Since then, the field of AI for mathematical discovery from Platonic data quickly blossomed and matured. Hundreds of papers emerged, with notable successes including the Ramanujan machine in finding identities for famous constants \cite{raayoni2021generating,raz2025euler}, the finding of new formulae in topology and representation theory \cite{davies2021advancing}, as well as the ML discovery of the still-open murmuration conjectures in number theory \cite{he2025murmurations}.

\section{Meta-Mathematics: }
Finally, I come to the perhaps most mysterious direction in which AI can galvanize mathematical research.
This is the world of Large Language Models (LLMs): chatGTP, Germini, Grok, Deepseek, etc.
It may seem peculiar that chatbots, famous for their hallucinations, could be of any use to mathematics. Indeed, though they have over the last couple of years rapidly out-performed and replaced humans in myriads of tasks
\footnote{Of course, it is a most welcome change that pointless
tasks such as form-filling and report-writing that dominate so much our day to day lives are now all done by LLMs.
It is remarkable that within two years we are beginning to end centuries of bureaucracy.}, LLMs struggle even in doing basic arithmetic. Language models rely entirely on statistical inferencing whereas mathematics requires {\em reasoning}.
Bottom-up is rigorous by construction, top-down's conjectures can be (dis-)proven rigorously. What can vague sentences and potential word-salads do for mathematical discovery? 

Curiously, in this direction, extraordinary advances are being made.
I will call this direction ``meta-'' since we could well be outside of mathematical reasoning, and depend solely on AI having linguistically seen (essentially all of) the mathematical literature.
From the first (and rather na\"{\i}ve) Word2Vec analyses of the ArXiv in order to uncover theories from language \cite{he2018hep} to the now much more sophisticated LLMMa2 \cite{touvron2023llama}, one should not underestimate the power of ``proof by vibe'' \cite{thurston2006proof}: an AI which has ``read'' all of contemporary mathematics and can be queried at will is a force to be reckoned with.
There is an increasing number of reports from mathematicians who, upon interacting with the likes of chatGTP5, have discovered useful literature of which they are not aware, or are instantly presented with ready-made computer code that turn out to be useful in testing ideas.
Meanwhile, LLMs can now attain the gold medal in the International Mathematics Olympiad by training on the language of previous solutions \cite{castelvecchi2025ai}. 
All of this is happening within the last couple of years.

In May of this year, some thirty of us were flown to the headquarters of EpochAI to help benchmark FrontierMath \cite{glazer2024frontiermath}, a project funded by OpenAI to see how various LLM models solve research-level mathematics, the sort of problems that one would work on and write papers about.
Accuracies of 10 - 30 \% have been achieved and the consensus was that oftentimes it truly felt like talking to an experienced colleague \cite{FT2025}.
The few hundred questions were divided into tiers, from undergraduate level (Tier 1) to publishable results (Tier 4), samples and performance statistics can be seen at \url{https://epoch.ai/frontiermath}.
The anticipated Tier 5 would involve open conjectures and problems, and the frightening thing is that we are not too far from this.

\section{Outlook}
In {\em Concrete Mathematics}, Graham, Knuth and Patashnik made the half-jocund remark  \cite{graham1994concrete} that  
``The ultimate goal of mathematics is to eliminate any need for intelligent thought''.
What they meant was that perhaps all methods of proof and derivation will be machanized in the infinite limit of mathematical research.
At the time, in the late 1990s, it was deemed outlandish such mechanization would ever be possible. Yet, the combination of human and AI in the last decade or so has given us {\em AI4Maths}, propelling us into a new era of fundamental theoretical discovery.
Like AlphaZero effectively \footnote{The ``effectively'' is important. Deepmind is not giving a deterministic solution or winning strategy, it is heuristically doing so by learning previous matches in AlphaGo and by playing against itself in AlphaZero.} solving Go, mathematics too can be gamified
\cite{shehper2024makes}.

Of course, we are still quite far from mechanized AI for mathematics which can follow an automated path of (1) find new problems from Platonic data and from existing literature; (2) find verifiable proof-paths using (vastly extended versions of) {\it Mathlib}; and (3) check the significance of the statement within the global and historical context.

Seeing that chatGTP has passed the Turing test \cite{biever2023chatgpt}, and that the AI community is working hard on the ARC test, a benchmark for AGI reasoning \cite{chollet2024arc}, a more stringent test was proposed for AI-assisted mathematical discovery \cite{he2024can}.
Named after Birch of the BSD Conjecture whose online talk inspired the participants of a conference at Cambridge in 2024 \cite{INI}, the Birch Test can be summarized as `` AIN@AI''.
An AI-guided mathematical discovery needs to be \underline{A}utomatic (in that human intervention is forbidden in choosing the data, the literature and the  algorithm), \underline{I}nterpretable (in that any statement must be meaningful to a human expert, e.g., having a deep NN that extrapolates an output without understanding the underlying mathematics is useless), and \underline{N}on-trivial (in that any conjecture raised must be significance to the scientific community).
While arguably only the Jones polynomial formula \cite{davies2021advancing}, the murmuration conjectures \cite{he2025murmurations} and finding unstable singularities for the Euler equation \cite{wang2025discovery} have come close to it, the Birch Test is so far unsurpassed \footnote{Or to some, perhaps even unsurpassable \cite{musaelian2025impossibility}.}.
Nevertheless, with the advent of {\em AlphaEvolve} \cite{novikov2025alphaevolve}, 
mathematical exploration and discovery are now at scale 
\cite{georgiev2025mathematical}. It should be noted here that the 67 problems selected for the study are of a particular nature - finding upper/lower bounds or numerical counter-examples - which are particular amenable to optimization and code-oriented searches.

Mathematics of the near future is clearly a joint effort between human expertise and AI's prowess in all three directions of bottom-up, top-down, and meta-mathematics.
Even in a hypothetical distant future where complete mechanization is achieved, where AI will (dis-)prove every major conjecture and propose new problems and continuously map out new areas of research, human mathematicians are still of great value. We will simply become priests to oracles, and interpret the results to the rest of humanity. Think of the philosophy departments in the world, centuries are spent in analyzing and critiquing Plato, or the literature departments, over Shakespeare. Perhaps one day in the far future, mathematics departments will consist of experts digesting the (Mathlib-verified) proofs that AI produces.
Still, I think a number of my colleagues agree with the exclamation \cite{HLF} that ``I don't care whether it is God, AI, or Terry Tao who finds the proof of the Riemann Hypothesis, I just want to know.''

\section*{Acknowledgments}
I am grateful to Ananyo Bhattacharya, Thomas Fink, Thomas Hodgkinson and Roman Rybiansky of the London Institute for Mathematical Sciences for their invaluable help and comments in the preparation of this Royal Institution lecture.
%%%%%%
\bibliography{refs}
\bibliographystyle{unsrt}
\end{document}